\documentclass[12pt,twoside,reqno]{amsart}
\usepackage{amscd,amssymb,amsmath}
\usepackage[arrow,matrix,curve]{xy}
\newtheorem{theorem}{Theorem}[section]
\newtheorem{proposition}[theorem]{Proposition}
\newtheorem{lemma}[theorem]{Lemma}
\newtheorem{corollary}[theorem]{Corollary}

\theoremstyle{definition}
\newtheorem{definition}[theorem]{Definition}
\newtheorem{example}[theorem]{Example}
\newtheorem{remark}[theorem]{Remark}

\newtheorem{construction}[theorem]{Construction}

\numberwithin{equation}{section}
\numberwithin{figure}{section}

\newcommand{\nc}{\newcommand}
\def\tag#1#2{\hbox to\textwidth{#1\hfil$\displaystyle #2$\hfil}}

\def\circ{\mathchoice%
 {\mathrel{\raise 1pt\hbox{$\scriptstyle\mathchar"020E$}}}
 {\mathrel{\raise 1pt\hbox{$\scriptstyle\mathchar"020E$}}}
 {\mathrel{\raise 1pt\hbox{$\scriptscriptstyle\mathchar"020E$}}}
 {}
}

\nc{\df}{{\it df}}
\nc{\dF}{{\it dF}}
\nc{\DF}{{\it DF}}
\nc{\ds}{{\it ds}}
\nc{\dvol}{{\it dvol}}
\nc{\grad}{{\it grad}}
\nc{\len}{{\it len}}
\nc{\PD}{{\it PD}}
\nc{\vol}{{\it vol}}
\nc{\area}{{\it area}}
\nc{\sys}{{\it sys}}
\nc{\fsys}{{\it fsys}}
\nc{\stsys}{{\it stsys}}
\nc{\confsys}{{\it confsys}}
\nc{\Hess}{{\it Hess}}
\nc{\parallelalf}{\parallel\!\alpha\!\parallel}
\nc{\parallelalfiminein}{\parallel\!I^{-1}(\alpha)\!\parallel}
\nc{\parallelalfpdr}{\parallel\!\PD_{\rr}(\alpha)\!\parallel}
\nc{\parallelpdreinsh}{\parallel\!\PD^{-1}_{\rr}(h)\!\parallel}
\nc{\parallelleer}{\parallel \; \parallel}
\nc{\parallelh}{\parallel\!h\!\parallel}
\nc{\parallelk}{\parallel\!k\!\parallel}
\nc{\parallelom}{\parallel\!\omega\!\parallel}
\nc{\parallelomx}{\parallel\!\omega_x\!\parallel}  
\nc{\parallelomi}{\parallel\!\omega_i\!\parallel}
\nc{\parallelpdralf}{\parallel\!\PD_{\rr}(\alpha)\!\parallel}
\nc{\parallelhfracnk}{\parallelh_{\frac{n}{k}}}
\nc{\strichleer}{| \:\: |}
\nc{\strichw}{|\omega|}
\nc{\strichwx}{|\omega_x|}
\nc{\strichalf}{|\alpha|}
\nc{\NN}{\mathbb N}
\nc{\RR}{\mathbb R}		
\nc{\R}{\mathbb R}
\nc{\rr}{\mbox{$\scriptstyle\mathbb R$}}
\nc{\nn}{\mbox{$\scriptstyle\mathbb N$}}
\nc{\zz}{\mbox{$\scriptstyle\mathbb Z$}}
\nc{\ZZ}{\mathbb Z}
\nc{\Z}{\mathbb Z}
\nc{\pisys}{{\it sys}\pi}

\def\C {{\mathbb C}} 
\def\AJ {{\mathcal A}} 
\def\ie {{\it i.e.\ }}
\def\cf {{\it cf.\ }} \def\eg {{\it e.g.\ }} \def\M {{\bf M}}  

\def\pvf{{- \hspace{-9.5pt}\mathcal{X}}} \def\genus{{s}}

\begin{document}

\author[V.~Bangert]{Victor Bangert} 
\address{ Mathematisches Institut,
Universit\"at Freiburg, Eckerstr.~1, 79104 Freiburg, Germany}
\email{bangert@mathematik.uni-freiburg.de}

\author[M.~Katz]{Mikhail Katz$^{*}$} \address{Department of
Mathematics and Statistics, Bar Ilan University, Ramat Gan 52900
Israel} \email{katzmik@math.biu.ac.il} \thanks{$^{*}$Supported by the
Israel Science Foundation (grant no.\ 620/00-10.0). Date: \today}

\title[Loewner inequality, harmonic 1-forms of constant norm]{An
optimal Loewner-type systolic inequality and harmonic one-forms of
constant norm}

\subjclass
{Primary 53C23;  
Secondary 57N65,  
52C07.		 
}

\keywords{Abel-Jacobi map, dual-extremal lattices, Loewner inequality,
minimal submanifold, Moser's volume-preserving flow, Riemannian
submersion, stable norm, systole}

\begin{abstract}
We present a new optimal systolic inequality for a closed Riemannian
manifold $X$, which generalizes a number of earlier inequalities,
including that of C.~Loewner.  We characterize the boundary case of
equality in terms of the geometry of the Abel-Jacobi map, $\AJ_X$, of
$X$.  For an extremal metric, the map $\AJ_X$ turns out to be a
Riemannian submersion with minimal fibers, onto a flat torus.  We
characterize the base of $\AJ_X$ in terms of an extremal problem for
Euclidean lattices, studied by A.-M.~Berg\'{e} and J.~Martinet.  Given
a closed manifold $X$ that admits a submersion $F$ to its Jacobi torus
$T^{b_1(X)}$, we construct all metrics on $X$ that realize equality in
our inequality.  While one can choose arbitrary metrics of fixed
volume on the fibers of $F$, the horizontal space is chosen using a
multi-parameter version of J.~Moser's method of constructing
volume-preserving flows.
\end{abstract}

\maketitle

\tableofcontents

\section{Introduction and statement of main theorems}  
\label{1}
We present a new optimal systolic inequality, which generalizes a
number of earlier inequalities.  These are the inequalities of
C.~Loewner \cite{Pu}, J.~Hebda \cite[Theorem A]{He}, certain results
of G. Paternain \cite{Pa}, and also inequality \cite[Corollary
2.3]{BK}.  Given a compact, oriented Riemannian manifold $X$ with
positive first Betti number $b>0$, our optimal inequality \eqref{11}
is a scale-invariant inequality providing a lower bound for the total
volume of $X$, in terms of the product of its conformal 1-systole and
its systole of codimension 1.  The definitions of the systolic
invariants involved appear in Section~\ref{222}.

In \cite[Corollary 2.3]{BK}, the authors generalize J.~Hebda's
inequality~\eqref{210} to arbitrary Betti number.  In the present
work, we strengthen the inequality by replacing the stable 1-systole
by its conformal analogue.

\begin{theorem}
\label{22}
Every compact, orientable $n$-manifold $X$ with positive first Betti
number $b=b_1(X,\RR)$ satisfies the inequality
\begin{equation}
\label{11}
\confsys_1(g) \sys_{n-1}(g)\le \gamma'_b \vol_n(g) ^{\frac {n-1} {n}}
\end{equation}
for every Riemannian metric $g$ on $X$.
\end{theorem}

The definition~\eqref{22'} of the Berg\'{e}-Martinet constant
$\gamma'_b$ and the associated dual-critical lattices, \cf
Definition~\ref{(2.one)}, implies that $\gamma'_b$ is precisely the
optimal constant in inequality \eqref{11} in the class of flat tori
$(T^b, g)$.  We also generalize Hebda's discussion of equality in
inequality~\eqref{210} as follows.

\begin{theorem}
\label{(2.2)}
Under the hypotheses of Theorem \ref{22}, one has equality in
\eqref{11} if and only if there exists a dual-critical lattice
$L\subset \R^b$ and a Riemannian submersion of $X$ onto the flat torus
$\R^b/L$, such that all fibers are connected minimal submanifolds of
$X$.
\end{theorem}

In particular, equality in \eqref{11} can only occur if $b_1(X,
\RR)\le n$.

An alternative way of stating the conclusion of Theorem~\ref{(2.2)} is
that we have equality in \eqref{11} if and only if the deck
transformations of the Jacobi variety \eqref{26} form a dual-critical
lattice, while the Abel-Jacobi map \eqref{ajm} is a Riemannian
submersion with connected minimal fibers.

\begin{example}
Let $(X,g)$ be a Riemannian 4-manifold with first Betti number
$b_1(X)=3$.  For instance, $X$ could be the product of a circle by the
nilmanifold given by the total space of the circle bundle over $T^2$
with nonzero Euler number.  Then inequality \eqref{11} takes the form
\[
\confsys_1(g)\sys_3(g)\leq \sqrt {\frac {3} {2} } \; \vol _4 (g)
^{\frac {3}{4}},
\]
with equality possible if and only if the deck transformations of the
Jacobi torus of $X$ form a cubic face-centered lattice, \cf
Example~\ref{22bis}.  The possible topology of $X$ was clarified in
\cite{NV}, \cf Remark~\ref{nagy} below.
\end{example}

Note that $X$ admits a Riemannian submersion with minimal fibers to a
flat $b$-dimensional torus if and only if all of its harmonic
one-forms have pointwise constant norm, \cf \cite{ES} and
Section~\ref{3}.  In Theorem~\ref{(7.3)} we characterize all the
metrics satisfying the boundary case of equality in \eqref{11}.  We
employ a multi-parameter generalisation of the method of
J.~Moser~\cite{Mo} of constructing volume-preserving flows to
determine all such metrics.

We study the relation between harmonic one-forms of constant norm and
Riemannian submersions to a flat torus with minimal fibers, \cf
Section \ref{3}.  Such submersions are special cases of foliations
with minimal leaves that have been studied in the theory of
foliations, \cf \cite{R}, \cite{Su}, \cite{KT}, and the survey by
V.~Sergiescu \cite{Se}.  However, the metrics in our context are much
more special.

In the class of flat tori $T^b$, the discussion of equality reduces to
an extremal problem for Euclidean lattices, studied by A.-M.~Berg\'{e}
and J.~Martinet \cite{BM}.  Namely, a flat torus $T^b$ realizes
equality if its deck transformations form a so-called dual-critical
lattice of rank $b$, \cf Definition~\ref{(2.one)}.  Thus, if $b=n$,
then Theorem \ref{(2.2)} implies that $(X,g)$ is a flat torus
$\R^n/L$, for a suitable dual-critical lattice $L\subset\R^n$.

If $\dim(X)=b_1(X)$ and all harmonic 1-forms have constant norm,
G. Paternain \cite[Corollary 2]{Pa} proved that $X$ is flat.  The
relation of the present work to \cite{Pa} is further clarified in
Remark~\ref{43bis}.

D.~Kotschick~\cite{Ko} defines the concept of {\em formal\/}
Riemannian metric and proves, in particular, that for such a metric
all harmonic forms (of arbitrary degree) have constant norm.
Moreover, if $n=b$ then $M$ is flat \cite[Theorem 7]{Ko}, while in
general one has $b\leq n$ \cite[Theorem~6]{Ko}.  Both of these
statements parallel our discussion of equality in \eqref{11}.

\begin{remark}
\label{nagy}
P. A. Nagy and C. Vernicos \cite{NV} study Riemannian manifolds
$(X,g)$ with all harmonic one-forms of constant norm in the case
$b_1(X)=\dim X-1$.  Assuming that $X$ is compact, orientable and
connected, they prove that $X$ is diffeomorphic to a two-step
nilmanifold.
\end{remark}

Systolic inequalities are metric-independent inequalities involving
minimal volumes of homology (or $1$-dimensional homotopy) classes in a
Riemannian manifold.  They have a certain similarity to the classical
isoperimetric inequality.  Over the past few years it has become clear
that systolic inequalities are systematically violated as soon as an
{\em unstable} $k$-systole with $k\geq 2$ is involved, \cf \cite{Ka1,
BabK, KS1, KS2, Bab, Ka2}.

On the other hand, a number of systolic inequalities, including those
involving stable and conformal systoles, may be found in M.~Gromov's
survey \cite{Gro2}, as well as in the more recent works \cite{BK, Ka3,
KKS}.  See the recent survey \cite[Figure~4.1]{CK} for a 2-D map of
systolic geometry, which places such results in mutual relation.

Systolic inequalities with {\it optimal constants} are rare.  In
addition to the ones by C.~Loewner \eqref{29} and J.~Hebda
\eqref{210}, there are the inequalities by P.M.~Pu \cite{Pu} for the
projective plane, by C.~Bavard \cite{Bav}, \cite{Sa} for the Klein
bottle, and by M.~Gromov \cite[pp.~259-260]{Gro3}, \cf
\cite[inequality (5.14)]{CK}.

\section{Norms on (co-)homology and definition of systoles}
\label{222}

In this section we assume that $X$ is a compact, oriented and
connected manifold, $\dim X=n$, with Riemannian metric $g$.  For $1\le
p\le\infty$ we define the $L^p$-norm $\parallelleer^*_p$ on
$H^k(X,\RR)$ and the dual norms $\parallelleer_p$ on the dual space
$H^k(X,\RR)^* = H_k(X,\RR)$.  Then $k$-systoles of $(X,g)$ will be
defined as the minima of such norms on non-zero elements in the
integer lattice $H_k(X,\ZZ)_{\RR}$ in $H_k(X,\RR)$.

The $L^p$-norms on $H^k(X,\RR)$, $k\in\{0,\ldots,n\}$, are quotient norms 
of the corresponding norms on closed forms. 
For $\omega\in\Omega^k X$ and $1\le p<\infty$ we set 
$$
\strichw_p = \left( \int_X \parallelomx^p \dvol_n(x)
\right)^{\frac{1}{p}},
$$
where $\| \omega_x \| =\max\{ \omega_x(v_1,\ldots,v_k)\mid |v_i|\le 1
\mbox{ for } 1\le i\le k \}$ is the point\-wise comass norm of
$\omega_x$.  In the cases $k=1$ and $k=n-1$ which are mainly of
interest in the present paper, this comass norm coincides with the
usual (Euclidean) norm $\strichwx$ induced by $g$.  For $p=\infty$, we
set
$$
\strichw_{\infty} = \max\limits_{x\in X}\parallelomx.
$$
If $\alpha\in H^k(X,\RR)$ and $1\le p\le\infty$, we define 
$$
\parallelalf^*_p = \inf \{\strichw_p \mid \omega\in\alpha\},
$$
where the expression $\omega\in\alpha$ means that $\omega$ is a smooth
closed $k$-form representing $\alpha$.  For $p=2$, Hodge theory
implies that every cohomology class $\alpha$ contains a unique
harmonic form $\omega$.  If $k\in\{1,n-1\}$, this form $\omega$ is
characterized by the equality $\strichw_2 = \parallelalf^*_2$.  In
particular, the norm $\parallelleer^*_2$ on $H^1(X,\RR)$ and
$H^{n-1}(X,\RR)$ is induced by the $L^2$ scalar product
$\langle,\rangle$, on harmonic forms.

The non-degenerate Kronecker pairing
$$
[,]: H_k(X,\RR) \times H^k(X,\RR)\rightarrow\RR
$$
induced by integration of closed forms over cycles, allows us to
define a norm $\parallelleer_p$ on $H_k(X,\RR)$ dual to the $L^p$-norm
$\parallelleer^*_p$ on $H^k(X,\RR)$.  Thus, given $h\in H_k(X,\RR)$,
we set
$$
\parallelh_p = \max\left\{ [h,\alpha] \left| \alpha\in H^k(X,\RR),
\parallelalf^*_p \le 1 \right. \right\} .
$$
Given a lattice $L$ in a finite-dimensional normed vector space
$(V,\parallelleer)$, we denote by
\begin{equation}
\label{(1.1)}
{\lambda_1(L) = \lambda_1(L,\parallelleer)}
\end{equation} 
the minimal norm of a non-zero lattice vector.  Then the {\it stable
$k$-systole} $\stsys_k(g)$ is defined as follows:
\begin{equation}
\label{(1.2)}
\stsys_k(g) = \lambda_1 \left( H_k^{\phantom{A}}(X,\ZZ)_{\RR},
\parallelleer_{\infty} \right).
\end{equation}
Here $H_k(X,\ZZ)_{\RR}$ denotes the lattice of integral classes in 
$H_k(X,\RR)$. 

The following geometric interpretation of the norm
$\parallelleer_{\infty}$ is well known and not obvious, \cf \cite[4.10
and 5.8]{Fe2} or \cite[4.18 and 4.35]{Gro3}.  Given an integer class
$h\in H_k(X, \ZZ)$, let $\vol_k(h)$ denote the infimum of the
Riemannian volumes of all integer $k$-cycles representing $h$.  Let
$h_\R$ denote the class corresponding to $h$ in $H_k(X,\R)$.
H. Federer proved the equality
\begin{equation}
\label{(1.3)}
\| h_\R \|_{\infty} = \lim\limits_{i\rightarrow\infty} \frac{1}{i}
\vol _k (ih),
\end{equation}
which is the origin of the term {\em stable norm\/} for $\parallelleer
_{\infty}$.

One also has the (homological) $k$-systoles, denoted $\sys_k$ and
defined by
\begin{equation}
\label{(1.4)}
\sys_k(g) = \min \left\{ \vol_k(h) \left|\; h\in H_k^{\phantom{A}} (X,
\ZZ) \setminus \{0\} \right. \right\} .
\end{equation}
If $H_k(X,\Z)$ is torsion-free, then $\stsys_k(g) \le \sys_k(g)$.  It
is a result due to H.~Federer \cite[5.10]{Fe2} that, in the assumption
that $X$ is orientable, one has
\begin{equation}
\label{(1.5)}
\parallelh_{\infty} = \vol_k(h) \mbox{ for all } h\in
H_{n-1}(X,\ZZ),
\end{equation}
and hence 
\begin{equation}
\label{(1.6)}
\stsys_{n-1}(g) = \sys_{n-1}(g).
\end{equation}
It is not difficult to see that $\sys_{n-1}(g)$ is the infimum of the
$(n-1)$-volumes of compact, orientable, non-separating hypersurfaces
in $X$.

For $p=\frac{n}{k}$, the norms $\parallelleer^*_p$ on $H^k(X,\RR)$ and
$\parallelleer_p$ on $H_k(X,\RR)$ are conformally invariant.  Hence
the quantity
\begin{equation}
\label{(1.7)}
\confsys_k(g) = \lambda_1\left( H_k(X,\ZZ)_{\RR}, \parallelleer
_{\frac{n}{k}} \right)
\end{equation}
is conformally invariant, and is called the {\it conformal}
$k$-systole of $(X,g)$.  If $h\in H_k(X,\RR)$ then the quantity $\| h
\|_{\frac{n}{k}}$ is the supremum of the stable norms of \nolinebreak
$h$ with respect to all metrics that are conformal to $g$ and have
unit volume, \cf \cite [7.4.A] {Gro1}.  In particular, one has
\begin{equation}
\label{28}
\stsys_k(g) \le \confsys_k(g) \vol_n(g)^{\frac{k}{n}}.
\end{equation}
However, the supremum of $\stsys_k(\tilde{g})$ over all unit volume
metrics $\tilde{g}$ conformal to $g$ may be smaller than
$\confsys_k(g)$.  Indeed, the usual systole and the conformal systole
obey dramatically different asymptotic bounds even in the case of
surfaces.  Thus, P. Buser and P.~Sarnak~\cite{BS} construct a sequence
of surfaces $\Sigma_\genus$ of genus $\genus\to\infty$, satisfying an
asymptotically {\em infinite\/} lower bound
\[
\confsys_1(\Sigma_\genus)^2\geq c\log \genus \mbox{ as } \genus \to
\infty,
\]
for a suitable constant $c>0$, while M. Gromov \cite{Gro2} brings the
asymptotically {\em vanishing\/} upper bound
\[
\frac{\sys_1(\Sigma_\genus)^2}{\area(\Sigma_\genus)}\leq C\frac {(\log
\genus) ^2}{\genus}.
\]
Now we can state the known systolic inequalities that are generalized
by Theorem~\ref{22}.

\subsubsection*
{\bf Loewner's inequality} (see \cite{Pu}) {\it Every riemannian
metric $g$ on the torus $T^2$ satisfies}
\begin{equation}
\label{29}
\sys_1(g)^2 \le \frac{2}{\sqrt{3}}\area(g).
\end{equation}
\quad A metric satisfying the case of equality in \eqref{29} is
homothetic to the quotient of $\C$ by the lattice spanned by the cube
roots of unity.

Note that in the case of the $2$-torus, the quantity $\sys_1(g)$ is
the length of the shortest noncontractible loop in $(T^2,g)$.

\subsubsection*
{\bf Hebda's inequality} (\cite[Theorem A]{He})  {\it Let $X$ be a
compact, orientable, $n$-dimensional manifold with first Betti number
$b_1(X,\RR)=1$.  Then every Riemannian metric $g$ on $X$ satisfies
\begin{equation}
\label{210}
\stsys_1(g) \sys_{n-1}(g) \le \vol_n(g)
\end{equation}
with equality if and only if $(X,g)$ admits a Riemannian submersion
with connected minimal fibers onto a circle.}

Note that inequality~\eqref{22} is a generalization of \eqref{29},
while \eqref{210} follows as a special case of inequality~\eqref{28}.
Indeed, the first Berg\'{e}-Martinet constant $\gamma'_1$ equals one,
while $\gamma'_2 = \frac{2}{\sqrt{3}}$, \cf
\cite[Proposition~2.13]{BM} and Section \ref{333}.

\section{Dual-critical lattices of Berg\'e and Martinet}
\label{333}
Theorem~\ref{22} involves a constant $\gamma'_b$, called the
Berg\'e-Martinet constant.  It is defined by setting, in the notation
of formula \eqref{(1.1)},
\begin{equation}
\label{22'}
\gamma'_b = \sup\left\{ \lambda_1(L) \lambda_1(L^*)\left| L \subseteq
(\RR ^b, \strichleer) \right. \right\},
\end{equation}
where the supremum is extended over all lattices $L$ in $\RR^b$ with a
euclidean norm $\strichleer$ and where $L^*$ denotes the lattice dual
to $L$.  The supremum defining $\gamma'_b$ is attained, \cf \cite{BM}.

\begin{definition}
\label{(2.one)}
A lattice $L$ realizing the supremum in formula~\eqref{22'} is called
{\em dual-critical}.
\end{definition}

\begin{remark}
\label{bmc}
The constants $\gamma'_b$ and the dual-critical lattices in $\R^b$ are
explicitly known for $b\le 4$, \cf~\cite[Proposition 2.13]{BM}.  In
particular, we have $\gamma'_1=1$, $\gamma'_2= \frac{2}{\sqrt{3}}$.
\end{remark}

\begin{example}
\label{22bis}
In dimension 3, the value of the Berg\'e-Martinet constant, $\gamma
_3'= \sqrt{\frac{3}{2}}=1.2247$, is slightly below the Hermite
constant $\gamma_3= 2^{\frac{1}{3}}=1.2599$.  It is attained by the
face-centered cubic lattice, which is not isodual \cite[p.~31]{MH},
\cite[Proposition 2.13(iii)]{BM}, \cite{CS}.
\end{example}

In general, the following facts are known about the constants
$\gamma'_b$:
\[
\gamma'_b \le \gamma_b \le \frac{2}{3} b \mbox{ for all $b\ge 2$}
\]
and
\[
\frac{b}{2\pi e} (1+o(1)) \le \gamma'_b\le\frac{b}{\pi e} (1+o(1))
\mbox{ as } b\rightarrow\infty,
\]
\cf \cite[pp.~334 and 337]{LLS}.

\section{Jacobi variety and Abel-Jacobi map}

Note that the flat torus $\R^b/L$ in Theorem~\ref{(2.2)} is isometric
to the Jacobi variety
\begin{equation}
\label{26}
J_1(X)=H_1(X,\R)/H_1(X,\ZZ)_\R
\end{equation}
of $X$, and the Riemannian submersion is the Abel-Jacobi map
\begin{equation}
\label{ajm}
\AJ_X: X \to J_1(X),
\end{equation}
induced by the harmonic one-forms on $X$, originally introduced by
A.~Lichnerowicz~\cite{Li}, \cf \cite[4.21]{Gro3}.

More precisely, let $E$ be the space of harmonic $1$-forms on $X$,
with dual $E^*$ canonically identified with $H_1(X,\R)$.  By
integrating an integral harmonic $1$-form along paths from a basepoint
$x_0\in X$, we obtain a map to $\R/\Z=S^1$. In order to define a map
$X\to J_1(X)$ without choosing a basis for cohomology, we argue as
follows.  Let $x$ be a point in the universal cover $\widetilde{X}$ of
$X$.  Thus $x$ is represented by a point of $X$ together with a path
$c$ from $x_0$ to it.  By integrating along the path $c$, we obtain a
linear form, $h\to \int_c h$, on $E$.  We thus obtain a map
$\widetilde{X}\to E^* = H_1(X,\R)$, which, furthermore, descends to a
map
\begin{equation}
\overline{\AJ}_X: \overline{X}\to E^*,\;\; c\mapsto \left( h\mapsto
\int_c h \right),
\end{equation}
where $\overline{X}$ is the universal free abelian cover.  By passing
to quotients, this map descends to the Abel-Jacobi map \eqref{ajm}.
The Abel-Jacobi map is unique up to translations of the Jacobi torus.

\section{Summary of the proofs}
\label{four}
The proof of Theorem~\ref{22}, which will be completed in
Section~\ref{nine}, depends on two results.  First we prove that for
conjugate exponents $p\in[1,\infty]$ and $q$, the Poincar\'e duality
map
$$
\PD_{\RR} : (H^1(X,\RR), \parallelleer^*_p)\rightarrow
(H_{n-1}(X,\RR),\parallelleer_q)
$$
is an isometry, see Proposition~\ref{41}.  On the other hand, the
H\"{o}lder inequality implies a chain of inequalities for the norms
$\parallelleer^*_p$ on $H^1(X,\RR)$, see Proposition~\ref{43}, and,
dually, opposite inequalities for the norms $\parallelleer_p$ on
$H_1(X,\RR)$.  This allows us to reduce inequality~\eqref{11} to the
case $p=2$ where the norms are Euclidean, so that the
definition~\eqref{22'} of $\gamma'_b$ applies.  If one has equality in
\eqref{11}, and if $\alpha\in H^1(X,\ZZ)_{\RR}$ is a nonzero element
of minimal $L^2$-norm, then the chain of inequalities for
$\parallelalf^*_p$ reduces to equality for all $p$, and we conclude
that the harmonic representative of~$\alpha$ has constant norm.

Similarly, the $L^2$-dual of every nonzero $h\in H_1(X,\ZZ)_{\rr}$ of
minimal dual $L^2$-norm is represented by a harmonic one-form of
constant norm, as well.  Using a result of Berg\'{e} and Martinet
\cite{BM} on dual-critical lattices, we then prove that equality in
\eqref{11} implies that all harmonic one-forms have constant norm.
Finally, the Riemannian submersion is the Abel-Jacobi map $\AJ_X$ of
formula~\eqref{ajm}, see Section~\ref{3}.  

In the case $b_1(X)=1$, J.~Hebda \cite{He} remarked already that
equality in his inequality~\eqref{210} implies the existence of a
nonzero harmonic one-form with constant norm, and hence of a
Riemannian submersion with minimal fibers over a circle.

Although Riemannian submersions over flat tori with minimal fibers are
relatively simple objects, it is not immediately clear what choices
one has in constructing all of them.  In the case $b_1(X)=1$, J.~Hebda
\cite{He} describes the simplest class of examples, namely local
Riemannian products.  In Theorem~\ref{(7.3)}, we present a
construction that starts with a submersion $F:X\rightarrow T^m$ over a
flat torus, and characterizes all Riemannian metrics on $X$ for which
$F$ is a Riemannian submersion with minimal fibers.

More specifically, we can prescribe arbitrary metrics of fixed volume
on the fibers.  To determine the metric, one additionally has to
choose a suitable horizontal distribution.  This choice depends on the
choice of a linear map from $\R^m$ into the Lie algebra of vertical
vector fields on $X$ that preserve the volume elements of the fibers.
If the fiber dimension, $\dim X-m$, is greater than one, this Lie
algebra is infinite-dimensional.  In particular, in this case the
Riemannian metrics satisfying the case of equality in \eqref{11}
are much more flexible than local Riemannian products.

\begin{remark}
\label{hyp}
There are compact manifolds that admit a submersion to the circle, but
do not carry a Riemannian metric possessing a nontrivial local product
structure, \eg hyperbolic manifolds fibering over a circle. Thus, our
construction yields new topological types of examples even in the case
$b_1(X)=1$.
\end{remark}

\section{Harmonic one-forms of constant norm and flat tori}
\label{3}
Most of the material of this section appears in \cite
[pp.~127-128]{ES}.  We include the proofs of Lemma \ref{(3.1)} and
Proposition \ref{(3.2)} for convenience.  The proofs are mostly
straightforward calculations in local Riemannian geometry.

First we give a geometric characterisation of the existence of a
nonzero harmonic one-form of constant (pointwise) norm on a Riemannian
manifold $(X,g)$.

\begin{lemma}
\label{(3.1)}
If $\theta\not= 0$ is a harmonic one-form of constant norm, then the
leaves of the distribution $\ker\theta$ are minimal, and the vector
field $V_{\theta}$ that is $g$-dual to $\theta$ has geodesic flow
lines orthogonal to $\ker\theta$.  

Conversely, given a transversely oriented foliation of $X$ by minimal
hypersurfaces, such that the orthogonal foliation consists of
geodesics, there exists a unique harmonic one-form $\theta$ with
$|\theta|=1$, such that $\ker\theta$ is everywhere tangent to the
foliation by minimal hypersurfaces, and such that $\theta$ is positive
on the oriented normals to the leaves.
\end{lemma}

\begin{proof}
Since $d\theta=0$, we can locally find a primitive $f$ of $\theta$,
\ie~$\df=\theta$ and hence $\grad f=V_{\theta}$.  The condition
$d^*\theta=0$ translates into $\mbox{div}(\grad f)=\Delta f=0$.  Since
$|V_{\theta}|=|\theta|$ is constant, for all vector fields $W$ on $X$
we have
\[
0=g(\nabla_W V_{\theta}, V_{\theta}) = \Hess(f) (W,
V_{\theta})=g\left(W,\nabla_{V_{\theta}} {V_{\theta}}\right).
\]
Hence $\nabla_{V_{\theta}} {V_{\theta}}=0$, \ie~the flow lines of
$V_{\theta}$ are geodesics.  Since $V_{\theta}$ is normal to the
leaves, while $|V_{\theta}|$ is constant and $\nabla_{V_{\theta}}
{V_{\theta}}=0$, the condition $\mbox{div} (V_{\theta})=0$ is
equivalent to the fact that the mean curvature of the leaves of
$\ker(\theta)$ vanishes.  

Conversely, let $\theta$ be the one-form with $|\theta|=1$ that
defines the given foliation by minimal hypersurfaces and is positive
on the oriented normals.  We want to prove that $\theta$ is closed.
Since the orthogonal foliation is geodesic by assumption, we have
$\nabla_{V_{\theta}} {V_{\theta}}=0$.  If $h$ is a locally defined
regular function whose level sets are leaves of $\ker(\theta)$, then
$\grad h=\lambda V_{\theta}$ for some nowhere vanishing function
$\lambda$.  Then for all vector fields $W, Z$ we have
\[
\Hess(h)(W,Z)=W(\lambda) g (V_{\theta},Z) + \lambda g(\nabla_W
V_{\theta},Z).
\]
Hence the bilinear form
\[
(W,Z)\rightarrow g(\nabla_W V_{\theta}, Z)
\]
is symmetric if $W$ and $Z$ are orthogonal to $V_{\theta}$.
Since $\nabla_{V_{\theta}} {V_{\theta}}=0$ and
$0=W(|V_{\theta}|^2) = 2 g(\nabla_W V_{\theta}, V_{\theta})$,
we conclude that
\[
g(\nabla_W V_{\theta}, Z) = g(\nabla_Z V_{\theta}, W)
\]
for all vector fields $W, Z$.  This implies $d\theta=0$.  As above, we
see that the minimality of the leaves, together with the conditions
$|V_{\theta}|=1$ and $\nabla_{V_{\theta}} V_{\theta}=0$, imply
$\mbox{div} (V_{\theta})=0$, and hence $d^*\theta=0$.
\end{proof}

Recall that a closed one-form is called {\em integral\/} if its
integrals over arbitrary closed curves are integers.

\begin{proposition}
\label{(3.2)}
For Riemannian manifolds $(X,g)$ the following two properties are
equivalent:
\begin{itemize}
\item[(i)] There exists a flat torus $(T^m, \overline{g})$ and a
Riemannian submersion $F:(X,g)\rightarrow(T^m,\overline{g})$ with
minimal fibers.
\item[(ii)] The set of harmonic one-forms of constant norm contains an
$m$-di\-men\-sional vector space generated by integral one-forms.
\end{itemize}
\end{proposition}

\begin{proof}[Proof of the implication $(i)\Rightarrow (ii)$]
Since $\overline{g}$ is flat, every harmonic one-form
$\overline{\theta}$ on $(T^m, \overline{g})$ has constant norm.
Integral forms pull back to integral forms, and so it suffices to
prove that $F^* \overline{\theta}$ is a harmonic one-form of norm
$|F^*\overline{\theta}| = |\overline{\theta}|$.  This equality is true
since $F$ is a Riemannian submersion.  To show that
$F^*\overline{\theta}$ is harmonic, note that $F^*\overline{\theta}$
is closed and of constant norm.  This implies that the foliation
orthogonal to $\ker\left(F^*\overline{\theta}\right)$ is geodesic if
$\overline{\theta}\not= 0$, \cf~the proof of Lemma~\ref{(3.1)}.  By
Lemma~\ref{(3.1)}, it remains to show that the leaves of the foliation
$\ker\left(F^*\overline{\theta}\right)$ are minimal.  This is a
consequence of the minimality of the fibers of $F$, together with the
fact that the leaves of $\ker\left(F^*\overline{\theta}\right)$ are
geodesic in horizontal directions.  Namely, if $v\in\ker\left( F^*
\overline {\theta}|_x\right)$ is orthogonal to the fiber
$F_x=F^{-1}(F(x))$ through $x$ and if $c$ is the geodesic in $X$ with
$\dot{c}(0)=v$, then $\dot{c}(t)\in \ker\left(F^*\overline {\theta}
|_{c(t)}\right)$ for all $t$.
\end{proof}

\begin{proof}[Proof of the implication $(ii)\Rightarrow (i)$]
Let $\theta_1, \ldots, \theta_m$ be linearly independent integral
harmonic one-forms on $X$, such that all linear combinations
$\sum\limits^m_{i=1} r_i \theta_i$ with $r_i \in\R$ have constant
norm.  Note that this last property is equivalent to the constancy of
all scalar products $g^*(\theta_i, \theta_j)$ with respect to the dual
metric $g^*$.  Let $\left(\overline{g}_{ij}\right)$ denote the matrix
inverse to $\left( g^* (\theta_i^{\phantom{A}}, \theta_j) \right)$.
Since the $\theta_i$ are integral, there exist functions $F_i$:
$X\rightarrow\R/\ZZ$ such that $d F_i=\theta_i$.  Then the map
\[
F=(F_1, \ldots, F_m): X\rightarrow \R^m/\ZZ^m
\]
is a Riemannian submersion from $(X,g)$ to $\R^m/\ZZ^m$ endowed with
the constant metric $\overline{g}=\sum\overline{g}_{ij} d x_i\otimes d
x_j$.  It remains to prove that the fibers of $F$ are minimal.  We can
find a global orthonormal frame field $V_1, \ldots, V_m$ for the
horizontal bundle ${\mathcal H}= \ker(DF)^{\perp}\subseteq T X$ by
taking the $g$-dual of appropriate constant linear combinations of
$\theta_1, \ldots, \theta_m$.  Locally we can complete $V_1, \ldots,
V_m$ to an orthonormal frame field $V_1, \ldots, V_n$ for $T X$.
Since $V_1, \ldots, V_m$ are $g$-dual to harmonic one-forms, we know
that
\[
\sum\limits^n_{j=1} g\left( \nabla_{V_j} V_i, V_j \right) =\mbox{div}
V_i = 0
\]
for $1\le i\le m$, \cf~the proof of Lemma~\ref{(3.1)}.  On the other
hand, Lemma~\ref{(3.1)} implies $\nabla_{V_j}V_j=0$ for $1\le j\le m$,
and hence
\[
g\left( \nabla_{V_j} V_i, V_j \right) = V_j (g(V_i, V_j))=0
\]
for $1\le j\le m$, $1\le i\le n$.  The preceding equations imply that
\[
\sum\limits^n_{j=m+1} g\left( \nabla_{V_j} V_i, V_j \right)=0
\]
for $1\le i\le m$.  Since $V_1, \ldots, V_m$ are orthogonal to the
fibers of $F$, this says that the mean curvature vector of the fibers
of $F$ vanishes.
\end{proof}

\begin{remark}
\label{43bis}
In the case $m=dim(X)$, the statement of Proposition~\ref{(3.2)} also
follows from the arguments employed by G. Paternain to prove
\cite[Corollary 2]{Pa}.  These arguments are global and use the
solution of the E. Hopf conjecture by D. Burago and S. Ivanov.  They
do not apply in the case $m< \dim(X)$.
\end{remark}

In the corollary below, we show that a Riemannian submersion $F$ as in
Proposition~\ref{(3.2)}(i) is uniquely determined by the metric $g$ on
$X$ and by the induced map $F_*$:$H_1(X,\ZZ)\rightarrow H_1
(T^m,\ZZ)$.

\begin{corollary}
\label{(3.3)}
Suppose $I:(X_0, g_0)\rightarrow(X_1, g_1)$ is an isometry and for
$i=0,1$ we are given Riemannian submersions $F_i$ with minimal fibers
from $(X_i, g_i)$ to flat tori $\left(T^m, \overline{g}_i\right)$. If
there exists an automorphism $A$ of $H_1(T^m,\ZZ)$ such that
$(F_1\circ I)_* = A\circ (F_0)_*$ then there exists an isometry
$\overline{I}=\left(T^m, \overline{g}_0\right) \rightarrow \left(T^m,
\overline{g}_1\right)$ such that $\overline{I}\circ F_0 = F_1\circ
I$.
\end{corollary}

\begin{proof}
\begin{figure}
\[
\xymatrix{ X_0 \ar[rr]^{I} \ar[d]^{F_0} &&X_1 \ar[d]^{F_1}\\ \left (
T^m,\bar g_0 \right) \ar@{-->}[rr]^{\overline{I}} && \left( T^m,\bar
g_1 \right) }
\]
\caption{Construction of isometry $\overline{I}$}
\label{31}
\end{figure}
Let $L_i\subseteq\R^m$ be the lattice of deck transformations of
$\left(T^m, \overline{g}_i\right)$.  Then we can identify $A$ with a
linear map $J:\R^m\rightarrow\R^m$ such that $J(L_0) = L_1$.  If
\[
\overline{J}:\R^m/L_0\rightarrow\R^m/L_1
\]
denotes the diffeomorphism induced by $J$, then our assumption implies
that $(F_1\circ I)_* = \left(\overline{J}\circ F_0\right)_*$.  Hence,
if $\theta$ is a constant one-form on $\R^m/L_1$, then the closed
one-forms $(F_1 \circ I)^* \theta$ and $\left(\overline{J}\circ
F_0\right)^* \theta$ are cohomologous.  From the proof of Proposition
\ref{(3.2)}, we know that $(F_1\circ I)^* \theta$ and
$\left(\overline{J}\circ F_0\right)^* \theta$ are harmonic with
respect to $g_0$.  By the uniqueness of the harmonic representative of
a cohomology class, we conclude that
\[
(F_1\circ I)^* \overline{\theta} = \left(\overline{J}\circ
F_0\right)^* \overline{\theta}
\]
for every constant one-form on $\R^m/L_1$.  This implies the
existence of $v\in\R^m/L_1$ such that
\[
F_1\circ I = \overline{J} \circ F_0 +v.
\]
If we define $\overline{I}:\R^m/L_0\rightarrow\R^m/L_1$ by
$\overline{I}(y)=\overline{J}(y)+v$, then $\overline{I}$ satisfies
$F_1\circ I=\overline{I}\circ F_0$, as illustrated in Figure~\ref{31}.
The map $\overline{I}$ is a diffeomorphism and a (local) isometry
since $I$ is an isometry, and $F_0$ and $F_1$ are Riemannian
submersions.
\end{proof}

\section{Norm duality and the cup product}
\label{4}
In this section, we assume that $(X,g)$ is a compact, $n$-dimensional,
oriented, and connected Riemannian manifold.  For $k\in\{1,n-1\}$ we
consider the norms $\parallelleer^*_p$ on $H^k(X,\RR)$ and
$\parallelleer_p$ on $H_k(X,\RR)$ defined in Section~\ref{222}.

In Proposition \ref{41}, we prove that for every pair of conjugate
exponents $p\in[1,\infty]$ and $q$, the Poincar\'e duality map is an
isometry between $\left( H^1(X,\R), \|\;\|_p \right)$ and $\left
( H_{n-1}(X,\R), \|\;\|_q \right)$, respectively between $\left(
H^{n-1}(X,\R), \|\;\|_p \right)$ and $\left( H_1(X,\R), \|\;\|_q
\right)$.

The results of the present section will be used in Sections~\ref{nine}
and \ref{eleven} to prove the systolic inequality in Theorem \ref{22}
and to analyze the case of equality.  For completeness and future
reference, some of them are proved in greater generality than
necessary for these applications.

\begin{proposition} 
\label{41}
Let $p\in [1,\infty]$ and $q$ be conjugate exponents.  Then the normed
spaces $\left( H^1(X,\R),\|\;\|_p^*\right)$ and $\left( H^{n-1} (X,
\R), \|\;\|_q^*\right)$ are dual to each other with respect to the cup
product.
\end{proposition}

\begin{remark}
\label{42}
The claim of Proposition \ref{41} is equivalent to either of the
following two statements:

\begin{enumerate}
\item
The Poincar\'e duality map
\[
\PD_\R : \left( H^1(X,\R),\|\;\|^*_p \right) \to \left(
H_{n-1}^{\phantom{A}} (X,\R),\|\;\|_q \right)
\]
is an isometry.
\item
The Poincar\'e duality map
\[
\PD_\R : \left( H^{n-1}(X,\R),\|\;\|^*_q \right) \to \left( H_1
^{\phantom{A}} (X,\R), \|\;\|_p \right)
\]
is an isometry.
\end{enumerate}
\end{remark}

\begin{proof}[Proof of Proposition \ref{41}]
A key idea is the application of the Hodge star operator in a suitable
measurable context, \cf formula \eqref{63bis}.  Let $\alpha\in
H^1(X,\R)$, $\beta\in H^{n-1}(X,\R)$.  If $\omega$ (respectively,
$\pi$) is a closed form representing $\alpha$ (respectively, $\beta$),
then the cup product $\alpha\cup\beta$ satisfies
\[
\PD_\R(\alpha\cup\beta)=\int_X\omega\wedge \pi.
\]
Since one- and $(n-1)$-forms take values in the set of simple
covectors, we can use \cite[p.~32]{Fe1} and the H\"{o}lder inequality
to estimate
\[
\left|\int_X \omega\wedge \pi \right| \leq \int_X |\omega_x| \;
|\pi_x| d\vol_n(x) \leq |\omega|_p |\pi |_q .
\]
This proves
\begin{equation}
\label{(42)}
\PD_\R(\alpha \cup \beta) \leq \|\alpha\|^*_p \| \beta\|^*_q
\end{equation}
It remains to show that for every $\alpha\in H^1(X,\R)\setminus \{ 0
\} $, there exists $\beta\in H^{n-1}(X,\R) \setminus \{ 0 \}$ such
that 
\begin{equation}
\label{(43)}
\PD_\R(\alpha\cup \beta) = \|\alpha\|^*_p \|\beta\|^*_q,
\end{equation}
or, conversely, that for every $\beta\in H^{n-1}(X,\R)\setminus \{ 0
\}$, there exists $\alpha\in H^1(X,\R) \setminus \{ 0 \}$ such that
\eqref{(43)} holds.  Assume first that $2\leq p < \infty$.  It is
shown in \cite[Proof of Lemma 4.2]{Si} that there exists a one-form
$\omega$ in the $L^p$-closure of $\alpha$ such that $|\omega|_p=
\|\alpha\|_p^*$.  It follows from \cite[Theorem 5.1]{Ha} that the form
\begin{equation}
\label{63bis}
\pi=|\omega|^{p-2}(*\omega) 
\end{equation}
is weakly closed and determines a class $\beta \in H^{n-1}(X,\R)$ such
that $\|\beta\|^*_q= |\pi|_q$.  This implies $0\not= \alpha\cup \beta
= \|\alpha\|_p^* \|\beta\|_q ^*$.

If $1 < p < 2$ and hence $2 < q < \infty$, we can apply analogous
arguments starting with an arbitrary $\beta\in H^{n-1}(X,\R)\setminus
\{ 0 \}$.  We then obtain $\alpha \in H^1(X,\R)\setminus \{ 0 \}$ such
that \eqref{(43)} holds.  Finally we treat the case $p=\infty$ and
$p=1$.  The proof in the case $p=1$ and $q=\infty$ is completely
analogous.  According to Remark \ref{42}, it suffices to prove that
the map
\[
\PD_\R : \left( H^{n-1}(X, \R), \|\;\|^*_1\right) \to \left(
H_1(X,\R), \|\;\|_\infty \right)
\]
is an isometry.  From \eqref{(42)} we conclude that $\|\PD
_\R(\beta)\|_\infty \leq \|\beta\|^*_1$ for every $\beta \in
H^{n-1}(X,\R)$.  To prove the opposite inequality, we use some
elementary facts from geometric measure theory.  Recall that a closed
normal one-current is a linear functional $T: \Omega^1 X \to \R$ such
that $T(df)=0$ for every $f\in C^\infty(M,\R)$ and such that its mass
$\M(T)$, defined by
\[
\M(T) = \sup \{ T(\omega) \mid \omega\in \Omega^1(M), |\omega |_\infty
\leq 1 \},
\]
is finite.  Every closed normal one-current $T$ determines a homology
class $[T]\in H_1(X,\R)$ such that $\left[ [H], \alpha \right] =
T(\omega)$ whenever $\omega\in \alpha$.  Every closed $(n-1)$-form
$\pi$ defines a closed normal one-current $T_\pi$ by setting
\[
T_\pi(\omega)= \int_X \omega \wedge \pi .
\]
Such one-currents are called smooth.  One easily sees that
\[
\M(T_\pi)= |\pi |_1
\]
and 
\[
\PD_\R(\pi)= \pm [ T_\pi ]
\]
Functional analysis implies that for every $h \in H_1(X,\R)$, there
exists a closed normal one-current $T$ such that $[T] = h$ and $\M(T)=
\|h\|_\infty$, \cf~\cite[Section 3] {Fe2}.

We can smooth this $T$, \cf \cite[4.1.18]{Fe2}, and obtain a sequence
$\pi_i \in \PD_\R^{-1}(h)$ such that
\[
\lim_{i\to \infty} |\pi_i|_1 = \lim_{i\to \infty} \M(T_{\pi_i}) =
\M(T) = \|h \|_\infty.
\]
This implies $\| \PD_\R^{-1}(h)\|^*_1 \leq \|h \|_\infty$.
\end{proof}

\section{H\"{o}lder inequality in cohomology and case of equality}
\label{541}

The H\"{o}lder inequality implies inequalities between the $L^p$ norms
$\|\;\|_p$ for different values of $p$.  We will prove that the
harmonic representative of a cohomology class $\alpha\in H^1(X,\R)$
has constant norm if and only if for some (and hence for every)
$p\not= 2$, the inequality relating the $L^2$- and the $L^p$-norm is
an equality at $\alpha$.  For $p=\infty$, this follows from
\cite[Theorem C]{Pa}.

\begin{proposition}
\label{43}
For every $\alpha \in H^1(X,\R)$, the function
\[
p\in [1,\infty] \to \|\alpha\|^*_p\vol_n(g)^{-\frac{1}{p}}
\]
is weakly increasing.  For $p\in [1,\infty]\setminus \{ 2 \}$, the
equality $\|\alpha \|^*_p \vol_n(g)^{-\frac{1}{p}} = \|\alpha \|^*_2
\vol_n(g)^{-\frac{1}{2}}$ holds if and only if the harmonic
representative of $\alpha$ has constant norm.
\end{proposition}

\begin{proof}
The monotonicity follows from H\"{o}lder inequality.  Next, we assume
that the harmonic representative $\omega$ of $\alpha$ has constant
norm.  Due to the monotonicity, it suffices to prove that
\[
\| \alpha \|^*_\infty \leq \| \alpha \| ^* _2 \vol_n(g)^{-\frac{1}{2}}
\leq \|\alpha \|^* _1 \vol_n(g)^{-1}
\]
The first inequality follows directly from the constancy of $| \omega
_x |$.  To prove the second inequality, let $(\omega_i)_{i\in \NN}$ be
a sequence of representatives of $\alpha$ such that $\lim_{i\to
\infty} |\omega|_1 = \|\alpha\|^*_1$.  Since $\omega _i -\omega$ is
exact and $\omega$ is harmonic, we see that $\langle \omega_i -
\omega, \omega \rangle = 0$ Since $| \omega_x |$ is constant, we
conclude that
\[
|\omega |_2^2 = \langle \omega_i , \omega \rangle_2 \leq |\omega_i|_1
| \omega | _2 \vol_n(g)^{-\frac{1}{2}}
\]
and hence
\[
\| \alpha \|_1^* \vol_n(g)^{-\frac{1}{2}} \geq \| \alpha \|_2 ^* .
\]
Finally, we assume that $p \in [1,\infty ] \setminus \{ 2 \}$ and that
\begin{equation}
\label{(44)}
\| \alpha \|_p ^* \vol_n(g)^{-\frac{1}{p}} = \| \alpha \|_2 ^*
\vol_(n)^{-\frac{1}{2}} ,
\end{equation}
and prove that the harmonic representative $\omega$ of $\alpha$ has
constant norm. Due to the monotonicity, it suffices to treat the cases
$1< p < 2$ and $2< p <\infty$.  If $1< p <2$, the discussion of
equality in the H\"{o}lder inequality implies that $\omega$ has
constant norm.  If $2 < p < \infty$, then there exists a one-form
$\omega$ in the $L^p$-closure of $\alpha$ such that $ | \omega |_p =
\| \alpha \|_p ^*$, \cf~\cite[proof of Lemma 4.2]{Si}.  Now the
H\"{o}lder inequality and \eqref{(44)} imply
\[
\begin{aligned}
\displaystyle
|\omega | _2 \vol_n(g)^{-\frac{1}{2}} & \leq | \omega |_p \vol_n(g)
^{-\frac{1}{p}} 
\\&
= \| \alpha \|_p ^* \vol_n(g)^{-\frac{1}{p}} 
\\&
= \| \alpha \|_2 ^* \vol_n(g)^{-\frac{1}{2}} .
\end{aligned}
\]
Hence $\omega$ is the harmonic representative of $\alpha$, the first
inequality is an equality, and the norm of $\omega$ is constant.
\end{proof}

Let $\langle \; , \; \rangle_2 ^*$ denote the scalar product on
$H^1(X,\R)$, with norm $\|\; \|_2^*$.  Thus, if we identify
$H^1(X,\R)$ with the space of harmonic one-forms, then $\langle \; ,
\; \rangle_2 ^*$ corresponds to the $L^2$-scalar product.  We define
an isomorphism $ I:H_1(X,\R)\rightarrow H^1(X,\R) $ by
\begin{equation}
\label{(45)}
\langle I(h),\alpha\rangle_{L^2}^* =[h,\alpha] \quad \mbox{ for } h\in
H_1(X,\R), \alpha\in H^1(X,\R).
\end{equation}
Note that $ \|h\|_2 = \|I(h)\|_2^*$ for all $h\in H_1(X,\R)$.
Proposition \ref{43} implies the following corollary.

\begin{corollary}
\label{44}
For every $h\in H_1(X,\R)$, the function
\[
p\in [1,\infty] \to \| h \|_p\vol_n(g)^{\frac{1}{p}}
\]
is (weakly) decreasing.  For $p\in [1,\infty] \setminus \{2 \}$, the
equality $\|h\|_p \vol_n(g)^{\frac{1}{p}} = \| h \|_2 \vol_n(g)
^{\frac {1}{2}}$ holds if and only if the harmonic representative of
$I(h)$ has constant norm.
\end{corollary}

Similarly, we can combine Propositions \ref{41} and \ref{43} to obtain
the following corollary.

\begin{corollary}
\label{45}
For every $k\in H_{n-1}(X,\R)$, the function
\[
p\in [1,\infty] \to \| k \|_p\vol_n(g)^{\frac{1}{p}}
\]
is (weakly) decreasing.  For $p\in [1,\infty] \setminus \{2 \}$, the
equality $\|k \|_p \vol_n(g)^{\frac{1}{p}} = \| k \|_2
\vol_n(g)^{\frac{1}{2}}$ holds if and only if the harmonic
representative of $\PD_\R^{-1}(k)$ has constant norm.
\end{corollary}

\section{Proof of Theorem \ref{22}}
\label{nine}

We assume that $(X,g)$ is a compact, oriented, and connected
Riemannian manifold with $b=b_1(X,\R)\geq 1$.  We first show how the
results proved in Sections~\ref{4} and \ref{541} imply
Theorem~\ref{22}:
\begin{equation}
\label{(51)}
\confsys_1(g)\sys_{n-1}(g)\leq \gamma_b' \vol_n(g)^{\frac{n-1}{n}} ,
\end{equation}
where $n=\dim(X)$ and $\gamma_b'$ is the Berg\'e-Martinet constant.
Note that inequality~\eqref{(51)} is true also for disconnected
manifolds once we know \eqref{(51)} for connected ones.

\begin{proof}[Proof of Theorem \ref{22}]
The Euclidean vector spaces $\left( H_1(X,\R), \|\;\|_2 \right) $ and
$\left( H_{n-1}(X,\R), \|\;\|_2 \right)$ are dual to each other with
respect to the intersection pairing defined by
\[
h \cdot k = [ h, \PD_\R^{-1}(k) ].
\]
This follows from the fact that the map $\PD_\R: \left( H^1(X,\R),
\|\;\|_2^* \right) \to \left( H_{n-1}(X,\R), \|\;\|_2 \right) $ is an
isometry, \cf Remark \ref{42}.  Moreover, the pair of lattices
$H_1(X,\Z)_{\rr} \subset H_1(X,\R)$ and $H_{n-1}(X,\Z)_{\rr} \subset
H_{n-1}(X,\R)$ are dual to each other with respect to this pairing,
\cf \cite[Section 3]{BK}.  Hence, formula~\eqref{22'} for the
Berg\'e-Martinet constant $\gamma_b'$ implies that
\[
\lambda_1\left( H_1^{\phantom{A}}(X,\Z)_{\rr}, \|\;\|_2\right) \lambda
_1 \left( H _{n-1}^{\phantom{A}} (X,\Z)_{\rr}, \|\;\|_2\right) \leq
\gamma_b' ,
\]
where $b=b_1(X)$.  Invoking Corollaries \ref{44} and \ref{45}, we
obtain for all $p,p'\in [ 2,\infty]$:
\begin{equation}
\label{(52)}
\lambda_1\left( H_1^{\phantom{A}}(X,\Z)_{\rr}, \|\;\|_p\right)
\lambda_1\left( H_{n-1}^{\phantom{A}}(X,\Z)_{\rr}, \|\;\|_{p'}\right)
\leq \gamma_b' \vol_n(g) ^{1- \frac{1}{p}- \frac{1}{p'}} 
\end{equation}
Now we specify parameter values to $p=n$ and $p'=\infty$.  Then
inequality~\eqref{(52)} specifies to \eqref{(51)}, \cf definition of
$\confsys_1$ in formula (2.7), as well as (2.2) and (2.6).
\end{proof}

Note that the inequality $\strichalf_{L^2} \left(\vol_n(g) \right)
^{\frac{1}{2}} \ge \|PD_{\rr}\alpha \|_\infty$ for a class $\alpha\in
H^1 (X,\R)$ can also be proved more directly, as follows.  Let $\alpha
\in H^1(X,\ZZ)_{\R}$ be a nonzero element in the integer lattice, and
let $\omega\in\alpha$ be the harmonic 1-form.  Then there exists a map
$f: X\to S^1=\R/\ZZ$ such that $df=\omega$.  Using the Cauchy-Schwartz
inequality and the coarea formula, \cf \cite[3.2.11]{Fe1},
\cite[p.~267]{Ch2}, we obtain
\begin{equation}
\begin{aligned}
\displaystyle |\alpha|_{L^2} (vol_n(g))^{1/2} & =
|\omega|_{L^2}(vol_n(g))^{1/2} \\ & \geq \int_X |df| dvol_n \\ & =
\int_{S^1} vol_{n-1}(f^{-1}(t))dt \\ & \geq \|\PD(\alpha)\|_\infty,
\end{aligned}
\end{equation}
since the hypersurface $f^{-1}(t)$ is Poincar\'{e} dual to $\alpha$
for every regular value of $t$.

\section{Consequences of equality and criterion of dual-perfection}
\label{ten}
We now collect some consequences of equality in \eqref{(51)}.  From
Corollaries~\ref{44} and \ref{45}, we obtain Lemma \ref{51} below.
\begin{definition}
A vector $v\in (L, \|\;\|)$ is called {\em short\/} if $\| v \| =
\lambda _1 (L, \|\;\|)$.
\end{definition}

\begin{lemma}
\label{51}
Assume equality holds in \eqref{(51)}.  Then
\begin{enumerate}
\item
if $\alpha \in (H^1(X,\Z)_{\rr}, \|\;\|_2^*)$ is short, then the
harmonic representative of $\alpha$ has constant norm;
\item
if $h\in (H_1(X,\Z)_{\rr}, \|\;\|_2)$ is short, then the harmonic
representative of $I(h)\in H^1(X,\R)$ as in (8.2) has constant
norm.
\end{enumerate}
\end{lemma}

Given a pair of dual lattices $L$ and $L^*$ in Euclidean space
$(E,\|\;\|_2)$, consider the sets of short vectors
$$
\begin{array}{rcl}
S(L) &=& \left\{ \alpha\in L \left|\; \|\alpha
\|_{2_{\phantom{A}}}^{\phantom{A}} = \lambda_1(L)\right. \right\}, \\
S(L^*) &=& \left\{ \beta\in L^* \left| \; \| \beta \|_2^{\phantom{A}}
= \lambda_1(L^*) \right. \right\},
\end{array}
$$
so that $S(L)\cup S(L^*)\subset E$.  Let $Q(E)$ denote the vector
space of symmetric bilinear forms on $E$, and let
$\varphi:E\rightarrow Q(E)^*$ be the map defined by
$$
\varphi(\gamma)(q)=q(\gamma, \gamma)
$$
for $\gamma\in E$ and $q\in Q(E)$.  

\begin{definition} [\cite{BM}]
\label{dp}
A pair $(L, L^*)$ of dual lattices in Euclidean space is called {\em
dual-perfect} if the set $\varphi(S(L)\cup S(L^*))$ generates $Q(E)^*$
as a vector space.
\end{definition}

The following was proved in \cite[item 3.10]{BM}.

\begin{theorem}
\label{dcdp}
Every dual-critical pair is necessarily dual-perfect.
\end{theorem}

Note that a compact, oriented Riemannian manifold for which equality
holds in \eqref{(52)} is necessarily connected.

\begin{proposition}
\label{(5.4)}
Let $(X,g)$ be a compact, oriented Riemannian $n$-manifold with
positive first Betti number $b=b_1(X,\R)$.  Assume that
$$
\confsys_1(g) \sys_{n-1}(g) = \gamma'_b \vol_n(g)^{\frac{n-1}{n}}.
$$
Then every harmonic one-form on $(X,g)$ has constant norm.
\end{proposition}

\begin{proof}
By Hodge theory, the space $E$ of harmonic 1-forms on $X$ endowed with
the $L^2$ scalar product is canonically isometric to the space
$\left(H^1(X,\R), \|\;\|_2\right)$.  In $E$ we have the pair of dual
lattices $L$ and $L^*$ that correspond to $H^1(X,\ZZ)_{\rr}$ and
$I\left(H_1(X,\ZZ)_{\rr}\right)$.

In the terminology of \cite{BM}, our assumption says that the pair
$(L, L^*)$ is dual-critical.  By Theorem \ref{dcdp}, the pair $(L,
L^*)$ is dual-perfect in the sense of Definition~\ref{dp}.  In our
situation, we have a canonical map $p:X\rightarrow Q(E)$, $x\mapsto
p_x$, defined by
$$
p_x(\alpha, \beta)=\langle\alpha(x), \beta(x)\rangle
$$
for $\alpha, \beta\in E$.  By Lemma~\ref{51}, every $\gamma\in
S(L)\cup S(L^*)$ is a harmonic form of constant norm.  This means that
$\varphi(\gamma)\in Q(E)^*$ is constant on the image $p(X)\subseteq
Q(E)$ of $p$.  Since the elements $\varphi(\gamma)$, $\gamma\in
S(L)\cup S(L^*)$ generate $Q(E)^*$, every element of $Q(E)^*$ is
constant on $p(X)$.  In particular, if $\gamma\in E$ then
$$
\varphi(\gamma)(p_x) = p_x(\gamma, \gamma) = |\gamma(x)|^2
$$
does not depend on $x\in X$.
\end{proof}

\begin{remark}
In terms of matrix traces, the argument above can be paraphrased as
follows.  The vectors $s\in S(L)\cup S(L^*)$ give rise to rank~1
symmetric matrices $ss^t\in Sym(\R^b)$.  A quadratic form $q_A(x)$
defined by a symmetric matrix $A$ gives rise to a linear map $Q:
Sym(\R^b)\to \R$ by the formula
$$
Q(xx^t)= q_A(x)= x^t A x = {\it Tr}(x^t A x) = {\it Tr}(Axx^t).
$$
Therefore, the map $Q$ is uniquely determined by its values on the
elements $ss^t\in Sym(\R^b)$, which form a spanning set by
dual-perfection.  It follows that every harmonic 1-form on $X$ has
constant norm.
\end{remark}

\begin{remark}
The following question was asked by A. Katok in 2002.  Does there
exist an example of a Riemannian metric with a pair of harmonic
one-forms of constant norm, whose pointwise scalar product is not
constant?  Such an example is, in fact, constructed in \cite{NV},
showing that the Berg\'e-Martinet criterion (Theorem~\ref{dcdp})
cannot be bypassed in the proof of Theorem~\ref{(2.2)}.
\end{remark}

\section{Proof of Theorem~\ref{(2.2)}}
\label{eleven}
Assume first that there exists a Riemannian submersion $F$ with
connected, minimal fibers of a compact, oriented, $n$-dimensional
Riemannian manifold $(X,g)$, with first Betti number $b\ge 1$, onto a
flat torus $\R^b/L$ where $L$ denotes a lattice in Euclidean space
$(\R^b, \strichleer)$.  We will show that in such a situation, we
necessarily have
\begin{equation}
\label{(61)}
\begin{array}{rcl}
\stsys_1(g) & = & \lambda_1(L),
\\
\confsys_1(g) & = & \lambda_1(L)\vol_n(g)^{-\frac{1}{n}}, 
\\
\sys_{n-1}(g) & = & \lambda_1(L^*)\vol_n(g).
\end{array}
\end{equation}
If $L$ is dual-critical, formulas~\eqref{(61)} imply that $\stsys_1(g)
\sys_{n-1}(g) = \gamma_b' \vol_n(g)$ and $\confsys_1(g) \sys_{n-1}(g)
= \gamma_b' \vol_n(g)^{\frac{n-1}{n}}$.  First note that the harmonic
one-forms on $X$ are precisely the pull-backs of the constant
one-forms on $\R^b/L$, and hence all of them have constant norm,
\cf~Proposition~\ref{(3.2)}.  Now Corollaries \ref{44} and \ref{45}
imply that, up to the factor $\vol_n(g)^{\frac{1}{p} - \frac{1}{2}}$,
the $L^p$-norms on $H_1(X,\R)$ and on $H_{n-1}(X,\R)$ coincide with
the $L^2$-norms $\|\;\|_2$.  In particular, we have
\begin{equation}
\label{(5.6)}
\stsys_1(g)=\lambda_1 \left( H_1^{\phantom{A}}(X,\ZZ)_{\rr},
\|\;\|_2^{\phantom{A}} \right) \vol_n (g) ^{\frac{1}{2}},
\end{equation}
and
\begin{equation}
\label{63}
\confsys_1(g)=\stsys_1(g)\vol_n(g)^{-\frac{1}{n}},
\end{equation}
and furthermore,
\begin{equation}
\label{(5.7)}
\sys_{n-1}(g)= \lambda_1 \left( H_{n-1}^{\phantom{A}}(X,\ZZ)_{\rr},
\|\;\|_2^{\phantom{A}} \right) \vol_n (g) ^{\frac{1}{2}}.
\end{equation}
By Remark~\ref{42}, the Poincar\'{e} duality map
\[
\PD_{\rr}: H^1(X,\R)\rightarrow H_{n-1}(X,\R)
\]
preserves the $L^2$-norms.  Since it maps $H^1(X,\ZZ)_{\rr}$ onto
$H_{n-1}(X,\ZZ)_{\rr}$, \cf \cite[Section 3]{BK}, equation
\eqref{(5.7)} is equivalent to
\begin{equation}
\label{(5.8)}
\sys_{n-1}(g) = \lambda_1 \left( H^1(X,\ZZ)_{\rr}, \|\;\|_2^* \right)
\vol_n (g) ^{\frac{1}{2}}.
\end{equation}
Since the fibers of the submersion $F$ are connected, the induced
homomorphism $\pi_1(X)\to \pi_1(\R^b/L)$ is onto, and hence we obtain
an epimorphism $H_1(X,\ZZ)\to H_1(\R^b/L,\ZZ)$.  Since $b=b_1(X,\R)$,
the epimorphism $F_* : H_1(X,\ZZ)_{\rr}\rightarrow H_1(\R^b/L,
\ZZ)_{\rr}$ is an isomorphism.  By duality, the map $F^* : H^1(\R^b/L,
\ZZ)_{\rr}\rightarrow H^1(X,\ZZ)_{\rr}$ is an isomorphism.  Now
$H^1(\R^b/L, \ZZ)_{\rr}$ can be identified with the dual lattice,
$L^*$.  An element $\theta\in L^*$ satisfies
\[
\|F^*\theta\|_2^*=|\theta|\vol_n(g)^{\frac{1}{2}}.
\]
Hence, equations~\eqref{(5.6)}, \eqref{63}, and \eqref{(5.8)} imply
our claim, equation~\eqref{(61)}.

Conversely, we assume that $\confsys_1(g) \sys_{n-1}(g) = \gamma'_b
\vol_n(g)^{\frac{n-1}{n}}$.  Then Proposition~\ref{(5.4)} implies that
every harmonic one-form on $X$ has constant norm.  Integrating
harmonic one-forms that form a basis for the integer lattice in
$H^1(X,\R)$, we obtain a Riemannian submersion $F$ with minimal fibers
from $(X,g)$ onto a flat torus $\R^b/L$, \cf~\eqref{(3.2)}.  Since we
obtain $F$ from a basis of the integer lattice $H^1(X,\ZZ)_{\rr}$, we
see that
\[
F^* : H^1(\R^b/L,\ZZ)_{\rr}\rightarrow H^1(X,\ZZ)_{\rr}
\]
is an isomorphism.  It follows that $F$ induces an epimorphism at the
level of fundamental groups, and hence the fibers of $F$ are
connected.  Note that the first part of the proof implies the
identities $\confsys_1(g) = \lambda_1(L)\vol_n(g)^{-\frac{1}{n}}$ and
$\sys_{n-1}(g) = \lambda_1(L^*)\vol_n(g)$.  Hence $\lambda_1(L)
\lambda_1(L^*) = \gamma'_b$, and therefore $L$ is dual-critical.

\section{Construction of all extremal metrics}
\label{seven}

In this section we will present a construction that starts with a
submersion $F:X\rightarrow T^m$ of a compact manifold $X$ to a flat
torus, and provides all Riemannian metrics on $X$ for which $F$ is a
Riemannian submersion with minimal fibers.  If $X$ admits a submersion
to $T^b$, where $b=b_1(X,\R)$, this construction provides all metrics
on $X$ for which all harmonic one-forms have constant norm,
\cf~Proposition~\ref{(3.2)}.

In this construction, we can prescribe arbitrary metrics of fixed
volume on the fibers of $F$.  This does not determine the metric on
$X$ uniquely.  We additionally have to define a horizontal
distribution $\mathcal{H}$ transverse to the vertical distribution
$\mathcal{V}=\ker(\dF)$.  The choice of $\mathcal{H}$ is restricted by
the requirement that the horizontal lift of every flow on the base has
to preserve the volume elements of the fibers.  The necessity of this
condition follows from the requirement that all fibers be minimal,
\cf~Lemma~12.3.  We start by providing a multiparameter
variant of J.~Moser's theorem \cite{Mo} on the existence of
volume-preserving diffeomorphisms.  This will allow us to find the
horizontal distributions, starting with an arbitrary family of
metrics on the fibers of constant volume.

We first explain the setting and fix some notation.  Let $F: X
\rightarrow Y$ denote a submersion between compact, connected and
oriented manifolds $X$ and $Y$.
\begin{definition}
Let $ \pvf _F$ denote the set of {\em projectable vector fields on}
$X$.  Namely, a vector field $V$ on $X$ is in $\pvf_F$ if there
exists a vector field $W$ on $Y$ such that $\DF\circ V=W\circ F$.
\end{definition}
Such a vector field $W$ on $Y$ will be called the {\it projection of}
$V\in \pvf_F$ and it will be denoted by $W=F_*(V)$.  In the
literature on foliations, projectable vector fields are also called
``basic'' or ``foliate''.  Note that $\pvf_F$ is a Lie
subalgebra of the Lie algebra of all vector fields on $X$.
Furthermore, $F_*$ is a homomorphism from the Lie algebra
$\pvf_F$ to the Lie algebra $\pvf(Y)$ of all vector
fields on $Y$.  In terms of the flow, $\Phi$, of $V\in\pvf_F$,
and the flow, $\Psi$, of $F_*(V)=W$, the relation $\DF\circ V=W\circ
F$ is equivalent to
$$
F\circ\Phi_t = \Psi_t\circ F \quad \mbox{for all $t\in\R$}.
$$
In particular, if $y\in Y$ then $\Phi_t$ maps the fiber $F^{-1}(y)$ to
the fiber $F^{-1}(\Psi_t (y))$.  We denote the diffeomorphism
$\Phi_t|F^{-1}(y)$ by $\Phi_{t,y}$.  

In addition, we assume that on every fiber $F^{-1}(y)$, we are given a
volume element $\alpha^y \in\Omega^{n-m}(F^{-1}(y))$, where $\dim X=n$
and $\dim Y=m$, that varies smoothly with $y\in Y$.  We say that the
flow $\Phi$ of $V\in\pvf_F$ {\em preserves}
$\alpha=(\alpha^y)_{y\in Y}$ if 
\[
\Phi^*_{t,y}\alpha^{\psi_t (y)}=\alpha^y ,
\]
for all $(t,y)\in\R\times Y$.  We denote by $\pvf (\alpha)$
the Lie subalgebra of all projectable vector fields whose flows
preserve $\alpha$.

\begin{proposition}
\label{(7.one)}
The homomorphism $F_* : \pvf
(\alpha)\rightarrow\pvf(Y)$ is an epimorphism of Lie algebras.
\end{proposition}

This multiparameter version of J.~Moser's theorem \cite{Mo} can be 
proved along the lines of the original proof in \cite{Mo}.

Now we specialize the preceding discussion to the case of a Riemannian
submersion $F:X\rightarrow Y$ between Riemannian manifolds $X$ and
$Y$, where $\dim X=n$, $\dim Y=m$.  Actually, the discussion is purely
local.  Every flow $\Psi$ on $Y$ has a unique {\it horizontal lift}
$\Phi$ to $X$, \ie~$\Phi$ is the flow on $X$ whose flow lines are the
horizontal lifts of the flow lines of $\Psi$.  Equivalently, if $\Psi$
is the flow of the vector field $W$ on $Y$, then $\Phi$ is the flow of
the {\it horizontal lift} $V$ of $W$ defined by requiring that $V \in
\pvf_F$, $F_*(V)=W$ and that $V$ is everywhere horizontal,
\ie~orthogonal to the vertical distribution ${\mathcal V}=\ker(\dF)$.
Obviously, the horizontal lift of any flow on $Y$ maps fibers of $F$
to fibers of $F$.

\begin{lemma}
\label{(6.3)}
\it Suppose $F:X\rightarrow Y$ is a Riemannian submersion between
oriented Riemannian manifolds. Then all the fibers of $F$ are minimal
submanifolds if and only if the horizontal lift of every flow on $Y$
maps the volume elements of the fibers of $F$ to each other.
\end{lemma}

\begin{proof}
This follows from \cite[Proposition 1]{R}.
\end{proof}

\begin{remark}
\label{(6.4)}
Actually, the ``if''-part in Lemma~\ref{(6.3)} is true under the
following weaker assumption. For every $y\in Y$ there are $m$ vector
fields on $Y$ that are linearly independent at $y$ and whose flows
have horizontal lifts that preserve the volume elements of the fibers
of $F$.
\end{remark}

Let us now assume that $F$ is a submersion of a compact, oriented and
connected $n$-manifold $X$ onto a flat torus $\R^m/L$ where $L$ is a
lattice in Euclidean space $\R^m$.  Moreover, we assume that for all
$y\in\R^m/L$, we are given a volume element $\alpha^y$ on $F^{-1}(y)$
that depends smoothly on $y$.

We let 
${\pvf}_c (\alpha)$ denote the set of all vector fields $V$ 
$\in {\pvf}(\alpha)$ such that 
$$
\dF(V)=(\dF_1(V), \ldots, \dF_m(V)) \text{ is a constant vector in }
\R^m.
$$
Obviously ${\pvf}_c (\alpha)$ is a linear subspace of ${\pvf}
(\alpha)$ and
\[
l_F : {\pvf}_c (\alpha) \rightarrow\R^m, \; l_F(V)=\dF(V),
\]
is a linear map.  Proposition~\ref{(7.one)} implies
\begin{equation}
\label{(7.1)}
l_F : {\pvf}_c (\alpha) \rightarrow\R^m \mbox{ is onto}.
\end{equation}
To define a horizontal distribution ${\mathcal H}$, we will choose a
linear right inverse $H:\R^m\rightarrow {\pvf}_c (\alpha)$ to
$l_F$ and set
$$
{\mathcal H}_x = \{H(v)|_x \mid v\in\R^m\}\subseteq TX_x
$$
for $x\in X$.  Note that we have much freedom in the choice of $H$ if
$n-m>1$.  The kernel of $l_F$ consists of those vertical vector fields
$V$ that preserve each $\alpha^y$, \ie~on every fiber $F^{-1}(y)$ one
has a vector field $V^y= V|F^{-1}(y)$ whose flow preserves $\alpha^y$
and that varies smoothly with $y\in T^m$.  In particular, $\ker(l_F)
\subseteq {\pvf}_c (\alpha)$ is an infinite-dimensional Lie
algebra if $n-m>1$.

We are now in a position to describe our construction of metrics for
which a given submersion onto a flat torus is Riemannian and with
minimal fibers.

\begin{construction}
\label{(7.2)}
Let $X$ be a compact, oriented and connected $n$-manifold and
$F:X\rightarrow T^m=\R^m/L$ a submersion onto a flat torus.
\begin{itemize}
\item[1)] Choose a smooth family $g^y$, $y\in T^m$, of Riemannian
metrics on the fibers $F^{-1}(y)$ of $F$, such that the total Riemannian
volume of $(F^{-1}(y), g^y)$ does not depend on $y\in T^m$.
\item[2)] Let $\alpha=(\alpha_y)_{y\in Y}$ denote the family of
Riemannian volume forms on the fibers of $F$.  Choose a linear right
inverse $H:\R^m\rightarrow {\pvf}_c (\alpha)$ to $l_F$ and define the
subbundle ${\mathcal H}\subseteq TX$ by
$$
{\mathcal H}_x = \{H(v)|_x \mid v\in\R^m\}
$$
for $x\in X$.
\item[3)] Now define a Riemannian metric $g$ on $X$ by requiring:
\begin{itemize}
\item[(a)] The subbundles ${\mathcal V}=\ker(\dF)$ and ${\mathcal H}$
are orthogonal to each other.
\item[(b)] For all $x\in X$ the isomorphism $\dF_x|{\mathcal H}_x :
{\mathcal H}_x\rightarrow\R^m$ is an isometry.
\item[(c)] For all $y\in T^m$ the inclusion $(F^{-1}(y),
g^y)\rightarrow(X,g)$ is isometric.
\end{itemize}
\end{itemize}
\end{construction}

\begin{theorem}
\label{(7.3)}
The Riemannian metrics $g$ on $X$ arising from
Construction~\ref{(7.2)} are precisely those for which $F:X\rightarrow
T^m=\R^m/L$ is a Riemannian submersion with minimal fibers.
\end{theorem}

\begin{proof} 
Suppose $g$ arises from the construction. Then properties (a) and (b)
imply that $F$ is a Riemannian submersion.  The vector fields
$V_i=H(e_i)$, $1\le i\le m$, form an orthonormal frame field for
${\mathcal H} = {\mathcal V}^{\bot}$.  Because of (c) and ${\pvf}_c
(\alpha) \subseteq {\pvf} (\alpha)$, their flows preserve the volume
elements induced by $g$ on the fibers.  This implies that the fibers
are minimal submanifolds of $(X,g)$, \cf~Remark~\ref{(6.4)}.
Conversely, suppose $g$ is a Riemannian metric on $X$ such that
$F:X\rightarrow\R^m/L$ is a Riemannian submersion with minimal fibers.
We describe how the metrics on the fibers and the map $H$ have to be
chosen so that our construction yields $g$.  For the metrics $g^y$,
$y\in\R^m/L$, on the fibers we must obviously take the metrics induced
by $g$.  If $v\in\R^m$ let $\overline{\theta}$ be the constant
one-form on $\R^m/L$ that is dual to $v$ with respect to the scalar
product on $\R^m$ and let $V$ be $g$-dual to $F^* \overline{\theta}$.
Now Lemma~\ref{(6.3)} implies that $V\in{\pvf}(\alpha)$, where
$\alpha=(\alpha_y)_{y\in Y}$ denotes the volume elements induced by
$g$ on the fibers.  Moreover, $l_F(V)=v$, \ie~the (linear) map that
maps $v\in\R^m$ to $V\in {\pvf}_c (\alpha)$ is a right inverse
to $l_F$ which we take to be $H$.  Then the horizontal distribution
${\mathcal H}$ induced by $H$ is $g$-orthogonal to ${\mathcal
V}=\ker(\dF)$ and our construction applied to these choices of $g^y,
y\in\R^m/L$, and $H$ yields the metric $g$.
\end{proof}

Corollary~\ref{(3.3)} shows to what extent our construction leads to
non-isometric metrics $g$, and how it depends on the submersion $F$
and the flat torus $\R^m/L$.

If $X$ admits a metric such that all harmonic one-forms have constant
norm, then there exists a submersion $F:X\rightarrow T^b$ for
$b=b_1(X,\R)$, \cf~Proposition~\ref{(3.2)}.  Starting from such
submersions $F$ and from arbitrary flat metrics on $T^b$, our
construction yields all metrics on $X$ for which all harmonic
one-forms have constant norm.

Similarly, $X$ can support metrics $g$ for which equality holds in
\eqref{11} only if $X$ admits a submersion onto $T^b$ for
$b=b_1(X,\R)$.  Starting from such submersions and from flat metrics
on $T^b$ coming from dual-critical lattices $L\subseteq\R^b$, the
construction yields all metrics on $X$ satisfying equality
in~\eqref{11}.

\section*{Acknowledgments}
We are grateful to D. Ebin, E. Leichtnam, and J.~Lafontaine for
discussions concerning the smooth dependence on parameters in Moser's
method, \cf Proposition~\ref{(7.one)}.  We thank E.~Kuwert for
providing references concerning closed forms minimizing the $L^p$-norm
in their cohomology class.

\bibliographystyle{amsalpha}

\end{document}